\documentclass[12pt,a4paper]{article}
\usepackage{amssymb}
\usepackage{amscd}
\usepackage{amsmath}
\begin{document}
\title{On isomorphism of simplicial complexes and their related rings}
\author{Rashid Zaare-Nahandi\\ {Institute for Advanced Studies in Basic
Sciences, Zanjan, Iran}\\ {E-mail: rashidzn@iasbs.ac.ir}}
\date{}
\maketitle
\begin{abstract}
In this paper, we provide a simple proof for the fact that two
simplicial complexes are isomorphic if and only if their associated
Stanley-Reisner rings, or their associated facet rings are
isomorphic as $K$-algebras. As a consequence, we show that two
graphs are isomorphic if and only if their associated edge rings are
isomorphic as $K$-algebras. Based on an explicit $K$-algebra
isomorphism of two Stanley-Reisner rings, or facet rings or edge
rings, we present a fast algorithm to find explicitly the
isomorphism of the associated simplicial complexes, or graphs.
\end{abstract}

\section*{Introduction}

Let $X$ be a finite nonempty set. A simplicial complex $\Delta$ on
$X$ is a set of subsets of $X$ such that, for any $x\in X$,
$\{x\}\in \Delta$, and if $E\in \Delta$ and $F\subseteq E$, then
$F\in \Delta$. A set in $\Delta$ is called a face and a maximal face
in $\Delta$ is called a facet.

Let $X=\{x_1,\ldots,x_n\}$, and $\Delta$ be a simplicial complex on
$X$. Let $R=K[x_1,\ldots,x_n]$ be the polynomial ring in $n$
indeterminates and with coefficients in a field $K$. Let $I(\Delta)$
be the ideal in $R$ generated by all square-free monomials
$x_{i_1}\ldots x_{i_s}$, provided that
$\{x_{i_1},\ldots,x_{i_s}\}\not\in\Delta$. The quotient ring
$R/I(\Delta)$ is called the Stanley-Reisner ring of the simplicial
complex $\Delta$.

It is easy to see that any quotient of a polynomial ring over an
ideal generated by square-free monomials of degree greater than 1,
is the Stanley-Reisner ring of a simplicial complex. A natural
question arises: If two Stanley-Reisner rings are isomorphic, are
their corresponding simplicial complexes isomorphic? In 1996, W.
Bruns and J. Gubeladze proved that if the isomorphism of the rings
is a $K$-algebra isomorphism, then the corresponding simplicial
complexes are isomorphic [1]. In this paper we provide a simple
proof for this result.

In 2002, S. Faridi in [2] defined the notion of facet ideal for a
simplicial complex which is a generalization of the notion of edge
ideal for a graph defined by R. Villarreal [6].

\paragraph{Definition 1.} Let $\Delta$ be a simplicial complex on the set
$X=\{x_1,\ldots,x_n\}$. Let $F(\Delta)$ be the ideal of
$R=K[x_1,\ldots,x_n]$ generated by all square-free monomials
$x_{i_1}\ldots x_{i_s}$, provided that $\{x_{i_1},\ldots,x_{i_s}\}$
is a facet in $\Delta$. The quotient ring $R/F(\Delta)$ is called
the facet ring of the simplicial complex $\Delta$.

\paragraph{Definition 2.} Let $G$ be a finite, simple and undirected graph with vertex
set $\{x_1,\ldots,x_n\}$. The ideal $E(G)$ of $R=K[x_1,\ldots,x_n]$
generated by all square-free monomials $x_{i}x_{j}$, provided that
$x_{i}$ is adjacent to $x_{j}$ in $G$, is called the edge ideal of
$G$. The quotient ring $R/E(G)$ is called the edge ring of the graph
$G$.

A question similar to the case of simplicial complexes can be stated
for facet rings and edge rings. In 1997, H. Hajiabolhassan and M. L.
Mehrabadi in [3] proved that two graphs are isomorphic if and only
if their corresponding edge rings are isomorphic as $K$-algebras. In
this paper we will prove the statement for the facet rings and
conclude it for the edge rings.

Finally, we will present a fast algorithm which admits a $K$-algebra
isomorphism of two Stanley-Reisner rings, or two facet rings, or two
edge rings as input and returns explicitly the isomorphism of
corresponding simplicial complexes or graphs, as output.

\section*{The isomorphism}
First we prove some lemmas.

\paragraph{Lemma 1.}
Let $I$ and $J$ be ideals in $R=K[x_1,\ldots, x_n]$ and
$S=K[y_1,\ldots,y_m]$, respectively, both generated by monomials of
degree greater than 1. Let $\phi:{R}/{I}\to{S}/{J}$ be a $K$-algebra
isomorphism. Let $f$ be a monomial in $R$ and denote the image of
$f$ in $R/I$ by $\bar{f}$. If $\bar{f}$ is a zero-divisor in
${R}/{I}$, then the constant term of $\phi(\bar{f})$ in ${S}/{J}$ is
zero.

\paragraph{Proof.} Let $\bar{g}$ be a nonzero monomial in
${R}/{I}$ and $\bar{f}\bar{g}=0$. Let $\phi(\bar{f})=\bar{h}=h+J$,
$\phi(\bar{g})=\bar{l}=l+J$. Let $\bar h=\bar h_0+\bar
h_1+\cdots+\bar h_r$ and $\bar l=\bar l_i+\bar l_{i+1}+\cdots+\bar
l_t$ be homogeneous decompositions of $\bar h$ and $\bar l$ with
$\bar l_i\neq 0$.  If $\bar h_0 \neq 0$, then $\bar h\bar l=\bar
h_0\bar l_i + \mbox{monomials of higher degree}$. But $\bar f\bar
g=0$ and then, $0=\phi(\bar f\bar g)=\bar h\bar l$. That is, $h.l\in
J$, and hence $l_i\in J$, which is a contradiction.
$\Box$\\

For any $i$, $1\leq i\leq n$, let $L(x_i)$ denote the set of
$y_j$'s such that $\bar y_j$ appears in the linear part of
$\phi(\bar x_i)$ with nonzero coefficient.

\paragraph{Lemma 2.} With the above notations, if $\phi:R/I\to S/J$
is a $K$-algebra isomorphism, then $m=n$ and $L(x_i)\neq
\varnothing$, for each $i$, $1\leq i\leq n$.

\paragraph{Proof.} The ideal $J$ is generated by monomials of degree
greater than one, therefore, $(S/J)_1$, the degree one homogeneous
component of $S/J$, is an $m$-dimensional $K$-vector space with
basis $\{\bar y_1,\ldots, \bar y_m\}$. The map $\phi$ is surjective,
thus the set $\{(\phi(\bar x_1))_1, \ldots, (\phi(\bar x_n))_1\}$
generates $(S/J)_1$ as a vector space over $K$. Therefore, $n\geq
m$. The map $\phi$ is isomorphism thus $\phi^{-1}$ is surjective too
and therefore, $m\geq n$. For the last claim, note that if for some
$i$, $1\leq i\leq n$, $L(x_i)=\varnothing$, then $(\phi(\bar
x_i))_1=0$ and the set $\{(\phi(\bar x_1))_1, \ldots, (\phi(\bar
x_n))_1\}$ can not
generate an $n$-dimensional vector space. $\Box$\\




Let $\Delta_1$ and $\Delta_2$ be two simplicial complexes on sets
$\{x_1,\ldots,x_n\}$ and $\{y_1,\ldots,y_m\}$, respectively. Let
$K[\Delta_1]=K[x_1,\ldots,x_n]/I(\Delta_1)$ and
$K[\Delta_2]=K[y_1,\ldots,y_m]/I(\Delta_2)$ be the Stanley-Reisner
rings associated to $\Delta_1$ and $\Delta_2$.

\paragraph{Theorem 1.} With the above notations, $\Delta_1$ and
$\Delta_2$ are isomorphic as simplicial complexes if and only if
$K[\Delta_1]$ and $K[\Delta_2]$ are isomorphic as $K$-algebras.

\paragraph{Proof.} It is obvious that if $\Delta_1$ and
$\Delta_2$ are isomorphic as simplicial complexes, then
$K[\Delta_1]$ and $K[\Delta_2]$ are isomorphic as $K$-algebras. For
the converse, assume that $\phi:K[\Delta_1]\to K[\Delta_2]$ is a
$K$-algebra isomorphism. By Lemma 2, $m=n$. An isomorphism as
$\phi$, can be uniquely determined by images of $\bar x_i$,
$i=1,\ldots,n$, and $\phi(\bar f(x_1,\ldots,x_n))=\phi(f(\bar
x_1,\ldots,\bar x_n))=f(\phi(\bar x_1),\ldots,\phi(\bar x_n))$ for
any polynomial $\bar f$ in $K[\Delta_1]$. Let $\phi(\bar
x_i)=f_i(\bar y_1,\ldots,\bar y_n)$, $i=1,\ldots,n$. By Lemma 1 and
Lemma 2, for each $i$, $f_i$ has no nonzero constant term and has
nonzero linear part. Suppose $\phi^{-1}$ is inverse of $\phi$ and
let $\phi^{-1}(\bar y_i)=g_i(\bar x_1,\ldots,\bar x_n)$,
$i=1,\ldots,n$. Let $f_{i1}$ and $g_{i1}$ denote the linear
homogeneous component of $f_i$ and $g_i$, respectively:
\begin{gather*}
f_{i1}=a_{i1}y_1 + \cdots + a_{in}y_n, \ \ \ \ i=1,\ldots,n \\
g_{j1}=b_{j1}x_1 + \cdots + b_{jn}x_n, \ \ \ \ j=1,\ldots,n.
\end{gather*}
The equalities $\phi\circ\phi^{-1}(y_i)=y_i$ and
$\phi^{-1}\circ\phi(x_i)=x_i$ for $i=1,\ldots, n$, imply that
$$
\left[
  \begin{array}{ccc}
    a_{11} & \ldots & a_{1n} \\
    a_{21} & \ldots & a_{2n} \\
    \vdots &        & \vdots \\
    a_{n1} & \ldots & a_{nn} \\
  \end{array}
\right] \left[
  \begin{array}{ccc}
    b_{11} & \ldots & b_{1n} \\
    b_{21} & \ldots & b_{2n} \\
    \vdots &        & \vdots \\
    b_{n1} & \ldots & b_{nn} \\
  \end{array}
\right] = \left[
  \begin{array}{ccc}
    b_{11} & \ldots & b_{1n} \\
    b_{21} & \ldots & b_{2n} \\
    \vdots &        & \vdots \\
    b_{n1} & \ldots & b_{nn} \\
  \end{array}
\right] \left[
  \begin{array}{ccc}
    a_{11} & \ldots & a_{1n} \\
    a_{21} & \ldots & a_{2n} \\
    \vdots &        & \vdots \\
    a_{n1} & \ldots & a_{nn} \\
  \end{array}
\right] = I
$$
where $I$ is the identity matrix of order $n$. Therefore, in the
following matrix, sum of entries of each row and each column is 1.
$$ M(\phi)=
\left[
  \begin{array}{cccc}
    a_{11}b_{11} & a_{21}b_{12} & \cdots & a_{n1}b_{1n} \\
    a_{12}b_{21} & a_{22}b_{22} & \cdots & a_{n2}b_{2n} \\
    \vdots & \vdots &  & \vdots \\
    a_{1n}b_{n1} & a_{2n}b_{n2} & \cdots & a_{nn}b_{nn} \\
  \end{array}
\right]
$$
 It is well known
that, there is a transversal of length $n$ with nonzero elements in
the matrix $M(\phi)$. See for instance [4]. By a transversal we mean
a sequence of entries of the matrix with no common columns or rows.
In other words, a transitive is a term in expansion of determinant
of the matrix. Let $1j_1, 2j_2,\ldots,nj_n$ be indices of a
transversal with nonzero elements in $M(\phi)$. By a change of
indices of $y_i$'s and permuting corresponding columns of $M(\phi)$,
suppose that the nonzero transversal is the main diagonal. Under
this assumption, $y_i\in L(x_i)$ and $x_i\in L^{-1}(y_i)$, for
$i=1,\ldots,n$, where $L^{-1}(y_i)$ is the set of variables with
nonzero coefficient in the linear part of $\phi^{-1}(y_i)$. For a
subset of $\{x_1,\ldots,x_n\}$ as $F$, which is not in $\Delta_1$,
there is a minimal set $E\subseteq F$, where $E\not\in\Delta_1$ and
any proper subset of $E$ belongs to $\Delta_1$. Let
$\{x_{i_1},\ldots,x_{i_r}\}$ be a minimal set not belonging to
$\Delta_1$. Then, $x_{i_1}\cdots x_{i_r}$ is a generator in
$I(\Delta_1)$, that is, $\bar x_{i_1}\cdots \bar x_{i_r}=0$ in
$K[\Delta_1]$ and $\phi(\bar x_{i_1}\cdots \bar x_{i_r})=0$ in
$K[\Delta_2]$. This means that $\phi(\bar x_{i_1})\cdots \phi(\bar
x_{i_r})\in I(\Delta_2)$. The ideal $I(\Delta_2)$ is a homogeneous
and monomial ideal and therefore each homogeneous component and each
monomial of $\phi(\bar x_{i_1})\cdots \phi(\bar x_{i_r})$ belongs to
$I(\Delta_2)$. Therefore, the product $f_{i_11}\cdots f_{i_r1}$ is
in $I(\Delta_2)$. In this case, there are two possibilities:
\begin{itemize}
\item
$y_{i_1}\cdots y_{i_r} \in I(\Delta_2)$, which means that
$\{y_{i_1}, \cdots, y_{i_r}\}\not\in \Delta_2$;
\item
$y_{i_1}\cdots y_{i_r}$ is canceled by another monomial in
$f_{i_11}\cdots f_{i_r1}$.
\end{itemize}
We prove that the second case is not possible. The condition in case
2, is equal to say that, for some $j$, $1\leq j\leq r$, $y_{i_j}$
appears with nonzero coefficient in more than one $f_{i_l1}$. Let
$y_{t_1}^{\alpha_1}\cdots y_{t_s}^{\alpha_s}$ be a smallest monomial
with lexicographic order in the set of monomials with nonzero
coefficient in the expansion of $f_{i_11}\cdots f_{i_r1}$ with
$\{t_1, \ldots, t_s\}\subseteq \{i_1,\ldots,i_r\}$ and $s$ is as
small as possible. This monomial appears once in the expansion and
so can not be canceled by another monomial. Therefore,
$y_{t_1}^{\alpha_1}\cdots y_{t_s}^{\alpha_s}\in I(\Delta_2)$ and
since $(\Delta_2)$ is a radical ideal, then $y_{t_1}\cdots
y_{t_s}\in I(\Delta_2)$ and so, $\bar y_{i_1}\cdots\bar y_{i_r}\in
I(\Delta_2)$.
Therefore, $\{x_{i_1},\ldots,x_{i_r}\}\not\in\Delta_1$ implies that
$\{y_{i_1},\ldots,y_{i_r}\}\not\in\Delta_2$. With a similar argument
for $\phi^{-1}$, $\{y_{i_1},\ldots,y_{i_r}\}\not\in\Delta_2$ implies
$\{x_{i_1},\ldots,x_{i_r}\}\not\in\Delta_1$, that is,
$\Delta_1\cong\Delta_2$.
$\Box$\\

Note that, for any simplicial complex $\Delta$, the ideal
$I(\Delta)$ has no monomial of degree one, but in the case of facet
ideals, zero dimensional facets correspond to degree one monomials
in $F(\Delta)$. To use Lemma 1 in the proof of the next theorem, we
will assume that there is not any zero dimensional facet in
simplicial complexes. This does not reduce the generality of the
theorem. Because, two simplicial complexes are isomorphic if and
only if they have the same number of zero dimensional facets and the
parts without any zero dimensional facet are isomorphic.

\paragraph{Theorem 2.} Any two simplicial complexes $\Delta_1$ and
$\Delta_2$ are isomorphic if and only if their corresponding facet
rings are isomorphic as $K$-algebras.

\paragraph{Proof.} The ``only if" part is obvious. To prove the ``if" part,
let
$$K[x_1,\ldots,x_n]/F(\Delta_1)\stackrel{\phi}{\to}K[y_1,\ldots,y_m]/F(\Delta_2)$$
be a $K$-algebra isomorphism between the facet rings corresponding
to $\Delta_1$ and $\Delta_2$. Similar to the proof of Theorem 1, it
follows that $m=n$. By an appropriate change of indices, we may
assume that $\{y_{1},\ldots,y_{n}\}$ is the set corresponding to the
main diagonal which we assumed to be a transversal with nonzero
elements in the matrix $M(\phi)$. Then, $y_{i}\in L(i)$ and $x_i\in
L^{-1}(y_{i})$, for $i=1,\ldots,n$. Let $\{x_{i_1},\ldots,x_{i_s}\}$
be a facet in $\Delta_1$. Then, $x_{i_1}\cdots x_{i_s}$ is in the
minimal generating set of $F(\Delta_1)$ and $\phi(\bar
x_{i_1})\cdots\phi(\bar x_{i_s})\in F(\Delta_2)$. The same argument
as the proof of Theorem 1, implies that $y_{i_1}\cdots y_{i_s}\in
F(\Delta_2)$. If $y_{i_1}\cdots y_{i_s}\in F(\Delta_2)$ is not in
the minimal generating set of $F(\Delta_2)$, then, without loss of
generality, we may assume that $y_{i_1}\cdots y_{i_s-t}\in
F(\Delta_2)$, for some $t$, $1\leq t\leq s-1$. This means that,
$\phi^{-1}(y_{i_1})\cdots\phi^{-1}(y_{i_s-t})\in F(\Delta_1)$ and
therefore, $x_{i_1}\cdots x_{i_s-t}\in F(\Delta_1)$, which is a
contradiction to minimality of $x_{i_1}\cdots x_{i_s}$ in
$F(\Delta_1)$. Therefore, $\{y_{i_1},\ldots,y_{i_s}\}$ is a facet in
$\Delta_2$, and this gives a bijection between $\Delta_1$ and
$\Delta_2$ as an isomorphism of simplicial complexes. $\Box$\\

A simple and undirected graph $G$ can be regarded as a simplicial
complex with facets $\{x_i,x_j\}$, where $x_i$ is adjacent to $x_j$
in $G$. With this interpretation, the edge ideal of $G$ is the same
as its facet ideal. Therefore, we have the following result.

\paragraph{Corollary.} Let $G_1$ and $G_2$ be two graphs. $G_1$ and $G_2$ are
isomorphic as graphs if and only if their edge rings are isomorphic
as $K$ algebras. $\Box$\\

\section*{The algorithm}
In this section, we present a fast algorithm to construct explicitly
an isomorphism between two simplicial complexes or two graphs, when
a $K$-algebra isomorphism of their associated rings is given.

\paragraph{Algorithm 1.} Let $R=K[x_1,\ldots,x_n]/I$ and $S=K[y_1,\ldots,y_n]/J$
be two $K$-algebras and $I$ and $J$ be ideals generated by some
square-free monomials of degree greater than one. Any $K$-algebra
homomorphism $\phi: R\to S$ can be uniquely determined by the images
of $\bar x_1, \ldots, \bar x_n$. In the following
algorithm, we assume that $R$ and $S$ are Stanley-Reisner rings of some simplicial complexes.\\
{\bf Input}: $K$-algebras $R$ and $S$, and a $K$-algebra
isomorphism $\phi: R\to S$, \\
{\bf Output}: A simplicial complex $\Delta_1$ associated to $R$ and
a simplicial complex $\Delta_2$ associated to $S$ as their
Stanley-Reisner rings and a bijection $\psi:\{x_1,\ldots,x_n\}\to
\{y_1,\ldots,y_n\}$ which determines an isomorphism of $\Delta_1$
and $\Delta_2$.
\begin{itemize}
\item[ ]{\it Step 1.}
Construct a simplicial complex $\Delta_1$ with the underlying set
$\{x_1,\ldots,x_n\}$ and faces $\{x_{i_1},\ldots, x_{i_r}\}$ where
$x_{i_1}\cdots x_{i_r}$ is not divided by any of generators of $I$.
Construct $\Delta_2$ on the set $\{y_1,\ldots,y_n\}$ and faces
$\{y_{i_1},\ldots, y_{i_s}\}$ where $y_{i_1}\cdots y_{i_s}\not\in J$
\item[ ]{\it Step 2.} Find the matrix  $M(\phi)$
\item[ ]{\it Step 3.} Find a transversal with nonzero elements in $M(\phi)$
\item[ ]{\it Step 4.} Use the transversal in step 3 to construct the map
$\psi:\Delta_1\to\Delta_2$.
\end{itemize}

The proof of Theorem 1, guaranties the correctness of the algorithm.
It is known that finding a nonzero transversal in a matrix which has
such a transversal, has a polynomial time algorithm [5], and so, the
above algorithm is polynomial time too.

\paragraph{Algorithm 2.} In Algorithm 1, we may consider the
$K$-algebras $R$ and $S$ as facet rings of simplicial complexes, or
if $I$ and $J$ are generated by square-free monomials of degree 2,
as edge rings of graphs. Then, in Step 1, we must construct a
simplicial complexes $\Delta_1$ and $\Delta_2$ such that $R$ and $S$
are their facet rings, respectively. Following steps 2, 3, and 4, we
finally obtain an isomorphism of simplicial complexes or graphs.

\paragraph{Acknowledgement.}
This work has been conducted during the Special Semester on
Gr\"obner Bases, February 1 - July 31, 2006, organized by Radon
Institute for Computational and Applied Mathematics (RICAM),
Austrian Academy of Sciences, and Research Institute for Symbolic
Computation (RISC) of Johannes Kepler University, Linz, Austria. The
author wishes to thank his colleagues at these institutes for their
warm hospitality.

\end{document}